{\bf Max Dehn,  Un Math\'ematicien Aux Pr\'eoccupations Universelles }  
                                    
  \vskip .3cm              
     {\bf  G. Burde,  W. Schwarz, J. Wolfart }
  \vskip .6cm   
trad. libre  Antoine D. Coste, CNRS.  

 \vskip 1.cm  
De 1921 \`a 1935 Max Dehn (1878-1952) contribua par des impulsions 
 d\'ecisives   \`a la premi\`ere floraison des math\'ematiques 
 \`a Francfort. Ses id\'ees compl\`etement nouvelles 
en th\'eorie des noeuds et en topologie
nourrirent le d\'eveloppement des math\'ematiques bien au del\`a 
des fronti\`eres de l'Allemagne. 
 En 1935 la terreur des nazis mit brutalement fin  \`a 
son activit\'e  \`a Francfort.

D\'epos\'e  et d\'echu de 
 sa citoyennet\'e,   
apr\`es une fuite p\'erilleuse   \`a travers
 la Norv\`ege, la Finlande, 
l' Union Sovi\'etique et le Japon,
 il trouva finalement refuge \`a 62 ans           
aux Etats Unis d'Am\'erique.    
Cependant, il ne put  r\'eussir \`a y obtenir un poste 
\`a la mesure de ses capacit\'es.

Le cinquanti\`eme anniversaire de sa mort a  \'et\'e l'occasion 
d'un regard r\'etrospectif.  
 \vskip 1.cm

En \'erudit complet, Max Dehn poss\'edait des  
 connaissances solides de la   philosophie grecque  
et de la math\'ematique.

  Dans le cadre de l'activit\'e  qu'il d\'eploya 
  en ses ann\'ees \`a  Francfort , 
ses contributions significatives 
couvrent trois domaines: les fondements de la g\'eom\'etrie,
 la th\'eorie des groupes et la topologie. Son nom 
reste attach\'e \`a de nombreux concepts novateurs 
sur les fondements de la  G\'eom\'etrie 
   et combinatoire des groupes: Le lemme et 
la chirurgie de Dehn, Dehn's filling, les torons (terme 
import\'e du langage des alpinistes francophones 
qui d\'ecrit les "twists" d'une corde )
  et les invariants de Dehn.

A 23 ans,  il d\'epose en 1901 un m\'emoire d'habilitation, 
par lequel il est le premier math\'ematicien \`a r\'esoudre 
un des c\'el\`ebres  probl\`emes de Hilbert, le troisi\`eme:   
                
" se donner deux t\'etra\^edres 
de m\^eme surface \`a la base 
 et de m\^eme hauteur qui ne se d\'ecomposent pas en 
t\'etra\^edres congruents........."   
\vskip .2cm  

L'\'elite du raisonnement se rassemblait au tournant du si\`ecle
dans l'\'ecole  de Hilbert  \`a  l'universit\'e de G\"ottingen. 
Max Dehn y   \'etudia. Il fut ensuite assistant  \`a Karlsruhe 
et \`a M\"unster  o\`u il passa son   habilitation.           . 
Il enseigna aussi \`a   Kiel et Breslau  (aujourd'hui Wroclaw, 
 Pologne). Ses ann\'ees de 1915 \`a 1918 sur le front francais, 
 jou\`erent un r\^ole dans sa situation en 1933. 
En 1921 il succ\'eda   \`a 
Ludwig Bieberbach  \`a Francfort. 
\vskip .5cm       

Pour d\'ecrire ce ``temps des cerises" (avant 1929  et 1933), 
 on cite souvent les 
``Souvenirs d'Apprentissage" 
d'Andr\'e  Weil, qui a s\'ejourn\'e longuement 
  \`a Francfort en 1926, puis 1952:
            
"Autour de  Max Dehn, Ernst Hellinger, Paul Epstein 
et Otto Szasz 
s'\'etaient rassembl\'es. Siegel \'etait  le dernier 
  arriv\'e.  D'eux, je ne peux parler sans avoir au 
coeur un sentiment  de reconnaissance".     

Weil compare Max Dehn  \`a  Socrate et 
 aborde le concept de sagesse: 
``Pour un tel homme la v\'erit\'e est indivisible 
 et la math\'ematique 
est seulement un des nombreux miroirs 
o\`u elle se  r\'efl\'echit,
peut \^etre avec plus de puret\'e qu'ailleurs(...).

Hellinger avait des dispositions semblables,
 quoiqu'un peu 
moins ardent. Certainement, il ne pouvait
 ext\'erioriser la m\^eme 
 autorit\'e morale que celle que seul 
Max Dehn d\'egageait par sa pr\'esence.                
Epstein et Szasz les soutenaient puissamment".  
\vskip 0.5cm 
Willy Hartner se souvient:``  (...) il aimait improviser 
et se laisser emplir de la profusion des pens\'ees 
qui affluaient  \`a son esprit. Auditeurs inexp\'eriment\'es, 
 il nous   \'etait un peu difficile 
d'appr\'ecier toute sa maitrise. D\'ecourag\'e,  
je lui demandai donc un court entretien,
 qui dura bien deux heures dans la cantine 
des professeurs, o\`u \`a cause de l'inflation 
calamiteuse, le caf\'e nous co\^uta un milliard de marks.
Je fus  le plus agr\'eablement surpris de ce que Dehn 
se montr\^at si disponible face \`a  mes interrogations. 
 
Le reste de la conversation roula sur l'art, la musique,
les langues classiques et modernes et aussi finalement 
sur la politique, la situation dans laquelle nous nous 
trouvions. Ce fut le d\'ebut d'une amiti\'e qui devait 
durer toute la vie, et a \'et\'e  \'eprouv\'ee en  une 
p\'eriode encore bien plus p\'enible.  " 
 \vskip 0.5cm 
En effet Willy Hartner, plus tard professeur d'histoire des 
sciences, mit \`a l'abri en 1938 des 
coll\`egues israelites, dont Max Dehn. Apr\`es la guerre il 
se pr\'eoccupa de ce que beaucoup de ceux qui avaient d\^u fuir
 \`a l'\'etranger   
pussent reprendre leur chaire, ou au moins recevoir leur 
pension ce qui, comme pour  Dehn, 
 fut fr\'equemment tr\`es difficile.

\vskip 1.cm 

Max Dehn,  Ernst Hellinger, Paul Epstein et Otto Szasz
 \'etaient d'origine israelite. En application de la loi nazie 
du 7 avril 1933 dite de "restauration de la carri\`ere  
 publique asserment\'ee" (Wiederherstellung des 
Berufsbeamtentums), Szasz perdit imm\'ediatement 
son poste de professeur et fut contraint \`a l'exil aux USA,
 Dehn,  Hellinger et Epstein,  anciens 
"Frontk\"ampfer" de la premi\`ere guerre mondiale, 
purent rester en poste. Cependant sous la pression de 
l'administration de l'universit\'e, Dehn et Hellinger 
d\'eprogramm\`erent leurs cours du semestre d'\'et\'e 1933. 
Siegel essaya de faire  en sorte qu'ils participassent 
  au moins au s\'eminaire, mais au cours de l'ann\'ee 
suivante l'atmosph\`ere ne fit que s'alourdir.   
 
Dans un article sur Dehn dans le Frankfurter Allgemeine Zeitung 
du 8 juillet 1952, W. Hartner \'ecrit:``   
Alors que le math\'ematicien Vahlen, dont Dehn avait 
d\'emoli le livre ``G\'eom\'etrie abstraite" dans une critique 
de 1905, \'etait nomm\'e en 1935 au minist\`ere de l'\'education, 
 Dehn comprit que ses jours \`a l'universit\'e de 
 Francfort \'etaient compt\'es. 
 Ce qu'il pr\'evoyait   survint imm\'ediatement: 
  sous le pr\'etexte de mesures d'\'economie, 
il fut mis \`a la retraite".  
   
En 1938, durant les pogromes de la nuit de cristal,
 Dehn fut arr\^et\'e, mais 
relach\'e car les prisons \'etaient pleines,
 et se mit \`a l'abri chez Willy Hartner.

A ce propos Siegel \'ecrit: " (...) Autrefois, ceux qui avaient 
cette conception de l'honn\^etet\'e    \'etaient en  minorit\'e 
et c' \'etait donc poss\'eder du courage que d'accueillir 
 un de ceux qui \'etaient pourchass\'es par les nazis et
leur pouvoir".

\vskip 1.cm  
 Transitant par Copenhague, Dehn r\'eussit \`a s'enfuir en 1939, 
 et prit un poste de professeur \`a l'universit\'e technique 
de Trondheim.  
L'invasion de la Norv\`ege obligea  Dehn 
\`a se cacher une deuxi\`eme fois et enfin,
 au d\'ebut de 1941 il 
entreprit, avec l'aide financi\`ere  de coll\`egues norv\'egiens 
une fuite p\'erilleuse en voiture vers la 
fronti\`ere su\'edoise
et un p\'enible voyage \`a travers la Finlande, 
l'Union Sovi\'etique et le Japon. 

   \vskip .5cm 
       
Ernst David  Hellinger(1883 Sil\'esie-1950 Chicago), 
 apr\`es avoir pass\'e 
six semaines au camp de concentration de Dachau, 
fut lib\'er\'e grace \`a l'argent que  sa soeur put  
envoyer depuis les  USA o\`u elle s' \'etait enfuie 
pr\'ec\'edemment et l'y rejoignit (NdT: apr\`es  
avoir rencontr\'e  Siegel \`a Francfort).     

\vskip .5cm          
                                              
Paul Epstein, qui \`a Strasbourg avait perdu 
son poste cons\'ecutivement \`a la premi\`ere guerre 
mondiale qu'il avait faite du  cot\'e allemand, 
et s'\'etait install\'e 
\`a Francfort,  se sentit trop vieux et trop 
malade pour s'enfuir   dans les ann\'ees 30. 
 Lorsqu'il recut en ao\^ut 1939 
une assignation pour interrogatoire 
 de la Gestapo, il s'empoisonna. 
Par les souvenirs de Viktor Klemperer, nous savons que 
cette fin tragique n'a pas  \'et\'e un  cas isol\'e.

\vskip 1.cm 
Dehn enseigna \`a Pocatello, Idaho de 1941 \`a 1942, 
au Illinois Institute of Technology   de 1942 \`a 1943, 
fut tuteur (!) au St John's College d'Annapolis 
de  1943 \`a 1944 et enfin professeur au 
Black Mountain College, Caroline du Nord, o\`u il 
dispensait des cours sur les math\'ematiques, la 
philosophie et les langues anciennes  
 (NdT:  peut \^etre sur  Epict\`ete ou sur les 
 fragments de Chrysippe,
\'edit\'es en 1905 par H. von Arnim \`a Stuttgart), 
 jusqu'\`a sa retraite 
en 1952. Sa femme Antonie  contribua \`a cette \'epoque 
 par des travaux d'arts d\'ecoratifs \`a leur situation 
mat\'erielle.  
  La pension de professeur 
fut donn\'ee \`a   son \'epouse 
le 1er avril 1950,  vers\'ee jusqu'\`a sa mort 
 soudaine le 27 juin 1952.   
W. Hartner critique  dans le Frankfurter Allgemeine Zeitung 
 les actes pris encore en avril 1952 par la "Wiesbadener 
 Fachbeh\"orde f\"ur Wiedergutmachung" du gouvernement 
r\'egional de  Hesse.   "R\'eparation" est 
 un mot  inappropri\'e et ind\'ecent. 
\vskip 1.cm 

Dehn eut trois remarquables \'etudiants: 
             
Ott Heinrich Keller(1906-1990), qui fut 
\`a partir de 1952 professeur ordinaire 
\`a Halle (alors RDA).   
\vskip 0.5cm     
  
Wilhelm Magnus(1907-1990), qui finit sa carri\`ere 
\`a New York. 
 \vskip 0.5cm  
Ruth Moufang (1905-1977), qui joua un r\^ole 
des plus importants 
pour le d\'eveloppement de l'universit\'e de Francfort. 
Docteur en 1930,
 recue \`a l'habilitation en 1936,       
 l'autorisation 
de d\'ecider seule du contenu de ses cours 
(``venia legendi"  )    
 lui fut refus\'ee par l'administration nazie. Elle quitta  
alors l' universit\'e pour un poste dans l'entreprise Krupp 
 mais retourna apr\`es la guerre au laboratoire
de math\'ematique (``Mathematische Seminar" ). Elle fut la 
premi\`ere femme nomm\'ee \`a un poste de professeur de 
math\'ematique en Allemagne, \`a l'Universit\'e  J.W. Goethe.  

         \vskip 1.cm  
{\ Bibliographie }
        \vskip 1.cm 
Andr\'e Weil    ``Souvenirs d'apprentissage", Birkhauser 1993.  
   
W.  Hartner , texte  \'edit\'e par la 
pr\'esidence de l'Universit\'e de Francfort.  
     
R. Moufang, W. Magnus, ``Max Dehn zum Ged\"achtnis",
Math. Annalen 127  (1954), 215-227. 
   
Oeuvres compl\`etes de C. L.  Siegel , vol 3, 462-474, 
conf\'erence pour le cinquantenaire de l'Universit\'e Goethe. 
         
John Stillwell, ``Max Dehn", History of Topology,
 ed. by I.M. James, 
Elsevier Science, 965-978, (1999);
``Max Dehn and geometry", Math. Semesterberichte 49,
 145-152(2002). 
  
J. W.  Dawson junior, ``Max Dehn, G\"odel, and the trans-Siberian" 
Int. Math. Nachrichten 189,1-13 (2002).
 
V. Klemperer,``Ich will Zeugnis 
ablegen bis zum letzten (Je veux rendre t\'emoignage 
jusqu'\`a la fin)"  
 journal, 1933-1941,
1942-1945  Aufbau Verlag, Berlin 1995.

\vskip 1.cm

r\'esum\'e     libre d'un article  (en allemand) des professeurs    
\vskip 1.cm    
 {\bf G. Burde } (auteur avec H. Zieschang 
de "Knots", \'editions  De Gruyter  , 1985) , 
                  
 {\bf W. Schwarz}  (ancien pr\'esident 
de la Soci\'et\'e math\'ematique d'Allemagne, 1986-1987, 
porte parole de la Conf\'erence des d\'epartements universitaires de 
math\'ematique, puis Directeur du 
d\'epartement de math\'ematique(Dekan )   
\`a la Johann Wolfgang Goethe Universit\"at zu Frankfurt am Main)   , 
               
 {\bf J. Wolfart } ( auteur de "Zahlen Theorie
  und Algebra" Vieweg Studium , 
Wiesbaden 1996, ancien directeur du d\'epartement de math\'ematique 
\`a la JWGUni-FFM, Directeur des \'etudes ), 
\vskip .5cm 

 {\bf paru dans 

"Forschung Frankfurt am Main, Das Wissenschaftsmagazin",
 pages 85-89,   en   avril 2002.}  

\vskip 1.cm
trad. Antoine D. Coste, CNRS KAF 194245,
 Sciences Physiques et Math\'ematique, membre de l'association 
``X-r\'esistance". 
\vskip .5cm
remerciements du traducteur \`a C.P. Korthals-Altes, Jean Lascoux, 
 G. Charpak, Alex Grosmann, Vincent Pasquier,
 A. Polyakov et  Daniel Altschuler
 pour des pr\'ecisions sur 
les torons, et  au Fachbereich Math JWG Uni pour son hospitalit\'e.      

\vskip 3.5cm 
 {\bf  R\'ef\'erences additionnelles } 
  
Walter Kempowski, oeuvres romanesques.   
      
G. Toulouse, L. Koch-Miramond, ``Les scientifiques et les droits 
humains" ed. de la maison des sciences de l'homme, Paris 2003.
                                                                         
G. Charpak, ``La vie \`a fil tendu", ed. O. Jacob, Paris. 
       
\vskip .3cm         
 e-archives en anglais: 
                           
 http://www-history.mcs.st-andrews.ac.uk/Mathematicians/Dehn.html  
                  
 http://www-history.mcs.st-andrews.ac.uk/Mathematicians/Hellinger.html

\end